\documentclass[a4paper,twoside,reqno]{amsart}

\usepackage[pagewise]{lineno}
\usepackage[utf8]{inputenc}

\usepackage{authblk}
\usepackage{hyperref}
\hypersetup{urlcolor=blue, citecolor=red}
\usepackage{amsmath,amsthm,amssymb,xcolor,bbm,mathrsfs}
\usepackage[english]{babel}
\usepackage{fancyhdr}
\newcommand{\norm}[1]{\|#1\|}

\newcommand{\n}{\hspace*{-6pt}}

\newcommand{\N}{\mathscr{N}}
\newcommand{\B}{\textbf{B}}

\DeclareMathOperator{\diag}{diag}

\DeclareMathOperator{\Span}{Span}

\makeatletter
\@namedef{subjclassname@1991}{2020 Mathematics Subject Classification}
\makeatother

\title{Leader-Follower dynamics}
\author{Hsin-Lun Li}

\date{}
\email{hsinlunl@asu.edu}

\theoremstyle{definition}
\newtheorem{theorem}{Theorem}

\newtheorem{lemma}[theorem]{Lemma}

\begin{document}

\allowdisplaybreaks

\thispagestyle{firstpage}
\maketitle
\begin{center}
    Hsin-Lun Li
    \centerline{$^1$National Sun Yat-sen University, Kaohsiung 804, Taiwan}
\end{center}
\medskip

\begin{abstract}
The original Leader-Follower model, proposed in \cite{zhao2018understanding}, categorizes agents with opinions in $[-1,1]$ into a follower group, a leader group with a positive target opinion in $[0,1]$, and a leader group with a negative target opinion in $[-1,0]$. Leaders maintain a constant attraction to their target, blending it with the average opinion of their group neighbors at each update. Followers, on the other hand, have a constant attraction to the average opinion of their leader group's opinion neighbors, also integrating it with their group neighbors' average opinion. This model was numerically studied in \cite{zhao2018understanding}. 

This paper extends the Leader-Follower model to include a social relationship, variable degrees over time, high-dimensional opinions, and a flexible number of leader groups. We theoretically investigate conditions for asymptotic stability or consensus, particularly in scenarios where a few leaders can dominate the entire population.

\end{abstract}

\section{Introduction}  
The Leader-Follower model contains the Hegselmann-Krause model in~\cite{mHK,mHK2} which involves two types of individuals: leaders and followers. The authors in \cite{zhao2018understanding} proposed a leader-follower model that partitions agents whose opinion is in $[-1,1]$ to a follower group, a leader group with a positive target opinion in $[0,1]$ and a leader group with a negative target opinion in $[-1,0]$. Individual $j$ is an opinion neighbor of individual $i$ if their opinion distance does not exceed the confidence threshold of individual $i.$ If all individuals share the same confidence threshold, two are opinion neighbors if their distance does not exceed that threshold. A leader's opinion depends on the opinion neighbors in its group and its group target, while a follower's opinion depends on all opinion neighbors. Define $[n]=\{1,2,\ldots,n\}$. Say $N$ agents including $N_1$ followers, $N_2$ positive target agents and $N_3$ negative target agents, set as $[N_1]$, $[N_1+N_2]-[N_1]$ and $[N]-[N_1+N_2]$.
The mechanism is as follows:
\begin{equation}\label{PN LF model}
\begin{array}{l}
\displaystyle x_i(t+1)=\frac{1-\alpha_i-\beta_i}{|N_i^F(t)|}\sum_{j\in N_i^F(t)}x_j(t)+ \frac{\alpha_i}{|N_i^P(t)|}\sum_{j\in N_i^P(t)}x_j(t) \vspace{2pt}\\
\displaystyle \hspace{55pt}+ \frac{\beta_i}{|N_i^N(t)|}\sum_{j\in N_i^N(t)}x_j(t),\quad i=1,\ldots,N_1,\vspace{2pt}\\
 \displaystyle   x_i(t+1)=\frac{(1-w_i)}{|N_i^P(t)|}\sum_{j\in N_i^P(t)}x_j(t)+ w_i d,\quad i=N_1+1,\ldots,N_1+N_2, \vspace{2pt}\\
 \displaystyle x_i(t+1)=\frac{1-z_i}{|N_i^N(t)|}\sum_{j\in N_i^N(t)}x_j(t)+ z_i g,\quad i=N_1+N_2+1,\ldots,N,

\end{array}
\end{equation}
where 
$$\begin{array}{rcl}
\displaystyle x_i(t) &\n=\n& \hbox{opinion of agent $i$ at time $t$},\vspace{2pt}\\
  \displaystyle   d &\n\in\n& [0,1]\ \hbox{is the positive target opinion},\vspace{2pt} \\
 \displaystyle    g &\in& [-1,0]\ \hbox{is the negative target opinion}, \vspace{2pt}\\
\displaystyle \epsilon_i &\n=\n& \hbox{confidence threshold of agent }i,\vspace{2pt}\\
\displaystyle N_i^F(t)&\n=\n& \{j\in [N_1]:\norm{x_i(t)-x_j(t)}\leq \epsilon_i\}, \vspace{2pt}\\
  N_i^P(t)&\n=\n& \{j\in [N_1+N_2]-[N_1]:\norm{x_i(t)-x_j(t)}\leq \epsilon_i\}, \vspace{2pt}\\
\displaystyle   N_i^N(t)&\n=\n& \{j\in [N]-[N_1+N_2]:\norm{x_i(t)-x_j(t)}\leq \epsilon_i\}, \vspace{2pt}\\
\displaystyle   \alpha_i &\n = \n& \hbox{degree to the average opinion of agent $i$'s} \\
&&\hbox{positive target neighbors}, \vspace{2pt}\\
\displaystyle   \beta_i &\n = \n& \hbox{degree to the average opinion of agent $i$'s}\\
&&\hbox{negative target neighbors}, \vspace{2pt}\\
   w_i &\n = \n& \hbox{degree to the positive target of agent }i, \vspace{2pt}\\
\displaystyle   z_i &\n = \n& \hbox{degree to the negative target of agent }i, \vspace{2pt}\\
   \alpha_i,\ &\beta_i,&\ w_i,\quad z_i \in [0,1].
\end{array}$$ 
The authors in \cite{zha2020opinion} pointed out that it can be an application in e-commerce. The Leader-Follower model we investigate now includes the following:
\begin{itemize}
    \item There is a social relationship.
    \item The degree toward the average opinion of a group can vary over time.
    \item The number of leader groups is decidable.
    \item Opinions can be high dimensional.
\end{itemize}
Let $\N_i^S(t)$ be the collection of all social and opinion neighbors of individual $i$ in set $S$ at time $t$ and $x_i(t)\in \mathbb{R}^d$ be the opinion of individual $i$ at time $t$ where $x_i(0)$ is a random variable.  The leader group $L$ model with target $g\in \mathbb{R}^d$ is given by:
\begin{equation}\label{leader group model}
    x_i(t+1)=\frac{\alpha_i(t)}{|\N_i^L(t)|}\sum_{j\in \N_i^L(t)}x_j(t)+(1-\alpha_i(t))g,\quad i\in L
\end{equation}
where $\alpha_i(t)\in [0,1]$ is a random variable indicating the degree of individual $i$ toward the average opinion of its group neighbors at time $t$ and $L$ is the collection of all leader group members. Say the leader groups are $L_1,\ldots,\ L_m$ with targets $g_1,\ldots,\ g_m.$ The follower group $F$ model is given by:
\begin{equation}\label{follower group model}
    x_i(t+1)=\frac{\big(1-\sum_{k=1}^m\beta_i^k(t)\big)}{|\N_i^F(t)|}\sum_{k\in \N_i^F(t)}x_k(t)+\sum_{k=1}^m\frac{\beta_i^k(t)}{|\N_i^{L_k}(t)|}\sum_{k\in \N_i^{L_k}(t)}x_k(t), \quad  i\in F
\end{equation}
where $\beta_i^k(t)\in [0,1]$ is a random variable indicating the degree toward the average opinion of its social and opinion neighbors in leader group $L_k$ at time $t$ and $F$ consists of all followers. $\beta_i^k(t)=0$ if $\N_i^{L_k}(t)=\emptyset.$ Observe that~\eqref{leader group model} and \eqref{follower group model} reduce to~\eqref{PN LF model} when there are two leader groups and a follower group. \eqref{leader group model} reduces to the synchronous Hegselmann-Krause model when $\alpha_i(t)=1$ for all $i\in L$ at all times. Similarly, \eqref{follower group model} reduces to the synchronous Hegselmann-Krause model when $\beta_i^k(t)=0$ for all $k\in [m]$ and $i\in F$ all the time. In~\eqref{leader group model} and~\eqref{follower group model}, a leader's opinion depends on the social and opinion neighbors in its group and its group target. In contrast, a follower's opinion depends on all social and opinion neighbors.

Interpreting in a graph, a vertex represents an individual and an edge symbolizes a relationship between two individuals. Saying
\begin{itemize}
    \item $G(t)=(V,E(t))$ is the social graph at time~t with vertex set and edge set $V$ and $E(t)$ and
    \item $\mathscr{G}(t)=(V,\mathscr{E}(t))$ is the opinion graph at time~t with vertex set and edge set $V$ and $\mathscr{E}(t).$
\end{itemize}  Edge \((i,j) \in E(t)\) if individual \(i\) is socially connected with individual \(j\), or if individual \(j\) is a social neighbor of individual \(i\). Similarly, edge \((i,j) \in \mathscr{E}(t)\) if individual \(i\) is opinion connected with individual \(j\), or if individual \(j\) is an opinion neighbor of individual \(i\). We can interpret a social relationship with an undirected social graph if edge \((i,j) \in E(t)\) implies \((j,i) \in E(t)\). For instance, if individual \(i\) is a relative of individual \(j\), then the reverse is also true. However, not all social relationships are reciprocal. For example, if individual \(i\) knows individual \(j\), it does not necessarily imply that individual \(j\) knows individual \(i\). In such cases, we use a directed social graph to represent the relationship. On the other hand, if all individuals share the same confidence threshold in an opinion relationship, we can interpret this opinion relationship with an undirected graph. A graph is $\delta$-\emph{trivial} if the opinion distance between any two vertices does not exceed $\delta.$ Denote $\B(a,r)$ as the open ball centered at $a$ with radius $r,$ i.e., $\B(a,r)=\{ x: \norm{x-a}<r\}.$ A \emph{profile} $G\cap \mathscr{G}$ is the intersection of the social and opinion graphs.

\section{Main results}
Since leader groups are independent, and similarly for the follower groups, we respectively investigate sufficient conditions for asymptotic stability or a consensus in the leader group $L$ and the follower group $F$ with $m$ leader groups. The sufficient condition in Theorem~\ref{Thm:consensus equal to target} is independent of social and opinion relationships. In fact, even a slight tendency toward the target by all leaders in $L$ guarantees a consensus equal to the target. 

\begin{theorem}\label{Thm:consensus equal to target}
     There is a consensus equal to the target in~\eqref{leader group model} when $$\liminf_{t\to\infty}\max_{i\in L}\alpha_i(t)<1.$$
\end{theorem}

The sufficient condition in Theorem~\ref{Thm:asymptotic stability on L} assumes an undirected social graph and an undirected opinion graph on $L$. Specifically, the synchronous Hegselmann-Krause model meets this condition, thus ensuring asymptotic stability. Asymptotic stability of the synchronous Hegselmann-Krause model illustrates finite time convergence.

\begin{theorem}\label{Thm:asymptotic stability on L}
    Assume that the social graph and opinion graph are undirected on $L$, the social graph becomes constant after some time, and $$\sum_{t\geq 0}(1/\min_{i\in L}\alpha_i(t)-1)<\infty.\ \hbox{Then, asymptotic stability holds in~\eqref{leader group model}.}$$
\end{theorem}

The sufficient condition in Theorem~\ref{Thm:consensus on F} specifies that the social graph and opinion graph can be directed, provided that all followers are socially connected with a leader in each leader group. This condition also identifies circumstances under which a few leaders can dominate the entire population.

\begin{theorem}\label{Thm:consensus on F}
    Assume that all followers have one social neighbor in each leader group, that $$ \big\{x_i(t),\ g_k\big\}_{i\in (\bigcup_{k=1}^m L_k)\cup F,\ k\in [m]}\subset \B(g_j,\min_{i\in F}\epsilon_i)$$ for some $j\in[m]$ and $t\geq 0$, that $\beta_i^k(t)=\beta_i^k$ for all $k\in [m]$ and $i\in F$ all the time, and that 
    $$\max_{i\in F}(1-\sum_{k=1}^m\beta_i^k)<1\quad \hbox{and}\quad \sup_{s\geq t}\big\{\max_{i\in L_k;k\in [m]}\alpha_i^k(s)\big\}<1.$$ Then, $$\lim_{t\to\infty}\max_{i\in L_k;k\in [m]}\norm{x_i(t)-g_k}=0\quad \hbox{and}\quad \lim_{t\to \infty}\max_{i\in F} \norm{x_i(t)- \frac{\sum_{k=1}^m\beta_i^k  g_k}{\sum_{j=1}^m\beta_i^j} }=0.$$
\end{theorem}

The sufficient condition in Theorem~\ref{Thm:asymptotic stability} assumes an undirected social graph and an undirected opinion graph on $F$, allowing the social graph and opinion graph on leader groups to be directed. Specifically, the synchronous Hegselmann-Krause model satisfies this condition.
\begin{theorem}\label{Thm:asymptotic stability}
    Assume that the social graph and opinion graph are undirected on $F$, the social graph becomes constant after some time, and $$\sum_{t\geq 0}\max_{i\in F; k\in [m]}\beta_i^k(t)<\infty.\ \hbox{Then, asymptotic stability holds in~\eqref{follower group model}.}$$ 
\end{theorem}

\section{The leader group model}
All leader groups are independent. Therefore, we investigate the behavior of a leader group. Let $y_i(t)=x_i(t)-g,$ \eqref{leader group model} becomes 
\begin{equation}\label{quasi HK model}
    y_i(t+1)=\frac{\alpha_i(t)}{|\N_i^L(t)|}\sum_{j\in \N_i^L(t)}y_j(t)
\end{equation}
It is clear that \eqref{quasi HK model} is the synchronous Hegselmann-Krause model when $\alpha_i(t)=1$ for all $i\in L$ all the time. Asymptotic stability holding in \eqref{leader group model} is equivalent to holding in \eqref{quasi HK model}. 

\begin{lemma}
    We derive $x_i(t)\to g$ if $\limsup_{t\to\infty}\alpha_i(t)=0$ for all $i\in L.$
\end{lemma}

\begin{proof}
    It follows from the triangle inequality that
\[
\norm{y_i(t+1)} \leq \alpha_i(t) \max_{j \in \N_i^L(t)} \norm{y_j(t)}.
\]
Taking \(\limsup\) on both sides, we get
\[
\limsup_{t \to \infty} \norm{y_i(t+1)} \leq \max_{i \in L} \norm{y_i(0)} \limsup_{t \to \infty} \alpha_i(t) = 0.
\]
This indicates that $y_i(t)\to 0$ as $t\to\infty.$ Hence, $x_i(t)\to g$ as $t\to\infty.$
\end{proof}

\begin{lemma}\label{lemma:supermartingale}
    Let $Z_t=\max_{i\in L}\norm{x_i(t)-g}.$ Then, $(Z_t)_{t\geq 0}$ is nonincreasing and 
    \begin{equation}\label{Eq: difference}
        Z_t-Z_{t+1}\geq \big(1-\max_{i\in L}\alpha_i(t)\big)Z_t.
    \end{equation}
\end{lemma}

\begin{proof}
    By the triangle inequality, we get $$Z_{t+1}=\max_{i\in L}\norm{y_i(t+1)}\leq\max_{i\in L}\alpha_i(t)\max_{i\in L}\norm{y_i(t)}=\max_{i\in L}\alpha_i(t)Z_t.$$ It turns out that $(Z_t)_{t\geq 0}$ is nonincreasing and $$Z_t-Z_{t+1}\geq \big(1-\max_{i\in L}\alpha_i(t)\big)Z_t.$$
\end{proof}

Next, we show circumstances in which social relationships and opinion relationships do not influence the achievement of consensus.

\begin{proof}[\bf Proof of Theorem~\ref{Thm:consensus equal to target}]
    It follows from Lemma~\ref{lemma:supermartingale} that $(Z_t)_{t\geq 0}$ is a nonnegative supermartingale. Via the martingale convergence theorem, $Z_t$ converges to some random variable $Z_{\infty}$ with finite expectation as $t\to\infty.$ Letting $\alpha_t=\max_{i\in L}\alpha_i(t)$ and taking $\limsup$ on~\eqref{Eq: difference}, we derive $$0=\limsup_{t\to\infty}(Z_t-Z_{t+1})\geq Z_\infty\limsup_{t\to\infty}(1-\alpha_t)=(1-\liminf_{t\to\infty}\alpha_t)Z_\infty.$$ This implies $Z_\infty=0.$
\end{proof}

\begin{lemma}\label{lemma:Wt}
    Assume that the social graph is undirected on $L,$ the opinion graph is undirected on $L$ with confidence threshold $\epsilon$, and $E(t)\subset E(t+1).$ Let $W_t=\sum_{i,j\in L}(\norm{x_i(t)-x_j(t)}^2\wedge \epsilon^2)\vee\epsilon^2\mathbbm{1}\{(i,j)\notin E(t)\}.$  Then, we derive
    \begin{equation}
    \begin{array}{rl}
         \n W_t-W_{t+1}\geq\n &\displaystyle 4\sum_{i\in L}\norm{x_i(t)-x_i(t+1)}^2-4|L|^2(1/\min_{i\in L}\alpha_i(t)-1)\vspace{2pt}\\
         &\displaystyle \times\max_{i\in L}\norm{x_i(0)-g}\bigg(\max_{i\in L}\norm{x_i(0)-g}\vee\max_{i,j\in L}\norm{x_i(0)-x_j(0)}\bigg) . 
    \end{array}  
    \end{equation}
\end{lemma}

\begin{proof}
Let $\N_i^L=\N_i^L(t),\ \alpha_i=\alpha_i(t),\ x_i=x_i(t),\ x_i^\star=x_i(t+1),\ y_i=y_i(t),\ y_i^\star=y_i(t+1),\ E=E(t)$ and $E^\star=E(t+1).$ It turns out that
    \begin{align*}
        W_t-W_{t+1}=&\sum_{i\in L}\bigg\{\sum_{j\in \N_i^L}(\norm{x_i-x_j}^2-\norm{x_i^\star-x_j^\star}^2\wedge\epsilon^2)\\
        &\hspace{1cm}+\sum_{j\in(\N_i^L)^c}\big[\epsilon^2-(\norm{x_i^\star-x_j^\star}^2\wedge\epsilon^2)\vee\epsilon^2\mathbbm{1}\{(i,j)\notin E^\star\}\big]\bigg\}\\
        &\hspace{-1.6cm}\geq \sum_{i\in L}\sum_{j\in \N_i^L}(\norm{x_i-x_j}^2-\norm{x_i^\star-x_j^\star}^2)=\sum_{i\in L}\sum_{j\in \N_i^L}(\norm{y_i-y_j}^2-\norm{y_i^\star-y_j^\star}^2)\\
        &\hspace{-1.6cm}=\sum_{i\in L}\sum_{j\in \N_i^L}\big(2<y_i-y_i^\star,y_i^\star-y_j>-2<y_i^\star-y_j,y_j-y_j^\star>\big)\\
        &\hspace{-1.6cm}=2\sum_{i\in L}|\N_i^L|<y_i-y_i^\star,y_i^\star>(1-1/\alpha_i)-2\sum_{i\in L}\sum_{j\in \N_i^L}<y_i^\star-y_i,y_j-y_j^\star>\\
        &\hspace{-1.1cm}-2\sum_{j\in L}\sum_{i\in \N_j^L}<y_i-y_j,y_j-y_j^\star>\\
        &\hspace{-1.6cm}=2\sum_{i\in L}|\N_i^L|<y_i-y_i^\star,y_i^\star>(1-1/\alpha_i)+2\sum_{i\in L}\norm{y_i-y_i^\star}^2\\
        &\hspace{-1.2cm}-2\sum_{i\in L}\sum_{j\in \N_i^L-\{i\}}<y_i^\star-y_i,y_j-y_j^\star>-2\sum_{j\in L}|\N_j^L|<y_j^\star/\alpha_j-y_j,y_j-y_j^\star>\\
        &\hspace{-1.6cm}\geq 2\sum_{i\in L}|\N_i^L|<y_i-y_i^\star,y_i^\star>(1-1/\alpha_i)+2\sum_{i\in L}\norm{y_i-y_i^\star}^2\\
        &\hspace{-1.1cm}-2\sum_{i\in L}\sum_{j\in \N_i^L-\{i\}}\norm{y_i^\star-y_i}\norm{y_j^\star-y_j}\\
        &\hspace{-1.1cm}-2\sum_{j\in L}|\N_j^L|(1/\alpha_j-1)<y_j^\star,y_j-y_j^\star>+2\sum_{j\in L}|\N_j^L|\norm{y_j-y_j^\star}^2\\
        &\hspace{-1.6cm}=-4\sum_{i\in L}|\N_i^L|(1/\alpha_i-1)<y_i^\star,y_i-y_i^\star>+2\sum_{i\in L}\norm{y_i-y_i^\star}^2\\
        &\hspace{-1.1cm}+\sum_{i\in L}\sum_{j\in \N_i^L-\{i\}}\big[(\norm{y_i^\star-y_i}-\norm{y_j^\star-y_j})^2-\norm{y_i^\star-y_i}^2-\norm{y_j^\star-y_j}^2\big]\\
        &\hspace{-1.1cm}+2\sum_{j\in L}|\N_j^L|\norm{y_j-y_j^\star}^2\\
        &\hspace{-1.6cm}\geq -2\sum_{i\in L}(|\N_i^L|-1)\norm{y_i^\star-y_i}^2-4\sum_{i\in L}|\N_i^L|(1/\alpha_i-1)<y_i^\star,y_i-y_i^\star>\\
        &\hspace{-1.1cm}+2\sum_{j\in L}|\N_j^L|\norm{y_j-y_j^\star}^2+2\sum_{i\in L}\norm{y_i-y_i^\star}^2\\
        &\hspace{-1.6cm}=4\sum_{i\in L}\norm{y_i-y_i^\star}^2-4\sum_{i\in L}|\N_i^L|(1/\alpha_i-1)<y_i^\star,y_i-y_i^\star>\\
        &\hspace{-1.6cm}\geq 4\sum_{i\in L}\norm{y_i-y_i^\star}^2-4|L|^2(1/\min_{i\in L}\alpha_i-1)\max_{i\in L}\norm{x_i(0)-g}\\
        &\hspace{3.3cm}\times\bigg(\max_{i\in L}\norm{x_i(0)-g}\vee\max_{i,j\in L}\norm{x_i(0)-x_j(0)}\bigg).
    \end{align*}
\end{proof}

By finiteness of the social graph, the social graph monotone after some time is equivalent to the social graph constant after some time. 

\begin{lemma}[Cheeger's Inequality \cite{beineke2004topics}]
\label{L7}
 Assume that~$G = (V, E)$ is an undirected graph with the Laplacian~$\mathscr{L}$. Define
 $$ i (G) = \min \bigg\{\frac{|\partial S|}{|S|} : S \subset V, 0 < |S| \leq \frac{|G|}{2} \bigg\} $$
 where~$\partial S = \{(u,v) \in E : u \in S, v \in S^c\}$. Then,
 $$ 2i (G) \geq \lambda_2 (\mathscr{L}) \geq \frac{i^2(G)}{2 \Delta (G)} \quad \hbox{where} \quad \Delta (G) = \ \hbox{maximum degree of~$G$}. $$
\end{lemma}

\begin{lemma}[\cite{mHK}]
\label{L10}
 Assume that~$Q$ is a real square matrix and that~$V$ is invertible such that the matrix~$VQ = \mathscr{L}$ is the Laplacian of some connected graph.
 Then, 0 is a simple eigenvalue of $Q'Q$ corresponding to the eigenvector~$\mathbbm{1} = (1, 1,\ldots, 1)'$.
 In particular, we have
 $$ \lambda_2 (Q'Q) = \min \{x'Q'Qx : \|x \| = 1 \ \hbox{and} \ x \perp \mathbbm{1} \}. $$
\end{lemma}

\begin{lemma}\label{lemma: delta-nontrivial component}
    Assume that the social graph and opinion graph are undirected on $L$. If some component $H$ of the profile $G(t)\cap\mathscr{G}(t)$ on $L$ is $\delta$-nontrivial, then
    \begin{align*}
        \sqrt{\sum_{i\in L}\norm{x_i(t)-x_i(t+1)}^2}&\\
        &\hspace{-1cm}\geq \sqrt{2}\delta\min_{i\in L}\alpha_i(t)/|L|^4-(1-\min_{i\in L}\alpha_i(t))\sqrt{|L|}\max_{i\in L}\norm{x_i(0)-g}.
    \end{align*}
\end{lemma}

\begin{proof}
    Letting \( V(H) \), the vertex set of \( H \), be \([h]\) and \(\mathbbm{1}=(1,\ldots,1)' \in \mathbb{R}^h\), express \(\mathbb{R}^h = W \oplus W^\perp\) for \( W = \Span(\{\mathbbm{1}\})\). For $y(t)=(y_1(t),\ldots,y_h(t))',$ write $$ y (t) = \left[c_1 \mathbbm{1} \,| \,c_2 \mathbbm{1} \,| \,\cdots \,| \,c_d \mathbbm{1} \right] +
            \left[\hat{c}_1 u^{(1)} \,| \,\hat{c}_2 u^{(2)} \,| \,\cdots \,| \,\hat{c}_d u^{(d)} \right] $$
    where $c_i$ and $\hat{c}_i$ are constants and $u^{(i)} \in \mathbbm{1}^\perp$ is a unit vector for all~$i \in [d]$. Observe that $$\norm{y_i(t)-y_j(t)}^2=\sum_{k\in [d]}\hat{c}_k^2(u^{(k)}_i-u^{(k)}_j)^2\leq 2\sum_{k\in [d]}\hat{c}_k^2\big((u^{(k)}_i)^2+(u^{(k)}_j)^2\big)\leq 2\sum_{k\in [d]}\hat{c}_k^2$$ for all $i,j\in [h].$
    Since $x_i-x_j=y_i-y_j$, $$\hbox{component}\ H\ \delta\hbox{-nontrivial implies}\ \sum_{k\in [d]}\hat{c}_k^2>\delta^2/2.$$ Letting $\alpha(t)=(\alpha_1(t),\ldots,\alpha_h(t))'$ and $B(t)=\diag(\alpha(t))A(t)$ for $A(t)\in \mathbb{R}^{h\times h}$ with $A_{i,j}(t)=\mathbbm{1}\{j\in \N_i^L(t)\}/|\N_i^L(t)|,$ we get $$y(t)-y(t+1)=(I-B(t))y(t)=\left[C(t)+F(t)\mathscr{L}(t)\right]y(t)$$ where $C(t)=I-\diag(\alpha(t))$, $F(t)=\diag(\alpha(t))\left(\diag((d_i)_{i=1}^h)+I\right)^{-1}$ with $d_i$ the degree of vertex $i$ in component $H$, and $\mathscr{L}(t)$ is the Laplacian of component $H.$ It follows from Lemmas~\ref{L7} and~\ref{L10} that
    $$\lambda_2(\mathscr{L})>\frac{(2/h)^2}{2h}=2/h^3,$$
    \begin{align*}
        \norm{F(t)\mathscr{L}(t)y(t)}^2&=\sum_{k\in [d]}\hat{c}_k^2\norm{F(t)\mathscr{L}(t)u^{(k)}}^2\geq \sum_{k\in [d]}\hat{c}_k^2\lambda_2\left(\mathscr{L}(t)F^2(t)\mathscr{L}(t)\right)\\
        &\geq (\delta^2/2) (\min_{i\in [h]}\alpha_i(t)/h)^2\lambda^2_2(\mathscr{L}(t))\geq  2\delta^2\min_{i\in [h]}\alpha_i^2(t)/h^8.
    \end{align*}
    On the other hand, we derive
    $$\norm{C(t)y(t)}\leq (1-\min_{i\in [h]}\alpha_i(t))\sqrt{h}\max_{i\in [h]}\norm{x_i(0)-g}.$$ 
    It follows from the triangle inequality that
    \begin{align*}
        \sqrt{\sum_{i\in L}\norm{x_i(t)-x_i(t+1)}^2}&\geq\sqrt{\sum_{i\in [h]}\norm{y_i(t)-y_i(t+1)}^2}=\norm{y(t)-y(t+1)}\\
        &\hspace{-2cm}=\norm{\left[F(t)\mathscr{L}(t)+C(t)\right]y(t)}\geq \norm{F(t)\mathscr{L}(t)y(t)}-\norm{C(t)y(t)}\\
        &\hspace{-2cm}\geq \sqrt{2}\delta\min_{i\in [h]}\alpha_i(t)/h^4-(1-\min_{i\in [h]}\alpha_i(t))\sqrt{h}\max_{i\in [h]}\norm{x_i(0)-g}\\
        &\hspace{-2cm}\geq \sqrt{2}\delta\min_{i\in L}\alpha_i(t)/|L|^4-(1-\min_{i\in L}\alpha_i(t))\sqrt{|L|}\max_{i\in L}\norm{x_i(0)-g}.
    \end{align*} 
\end{proof}


\begin{proof}[\bf Proof of Theorem~\ref{Thm:asymptotic stability on L}]
    We claim the following:
    \begin{enumerate}
        \item\label{C1} All components of profile $G\cap\mathscr{G}$ on~$L$ are $\delta$-trivial after some time for all $\delta>0.$\vspace{2pt}
        \item\label{C2} No components of profile $G\cap\mathscr{G}$ on~$L$ interact with each other after some time.
    \end{enumerate}
    Without loss of generality, we assume the social graph on $L$ remains constant over time, saying $G(t)|_L=G|_L=(L, E)$ for all $t\geq 0$. Observe that $$\sum_{t\geq 0}(1/\min_{i\in L}\alpha_i(t)-1)<\infty \implies \lim_{t\to\infty}\min_{i\in L}\alpha_i(t)=1\iff \lim_{t\to\infty}\alpha_i(t)=1\ \hbox{for all}\ i\in L.$$ Hence, we derive $$\sqrt{2}\delta\min_{i\in L}\alpha_i(t)/|L|^4\to\sqrt{2}\delta/|L|^4\ \hbox{and}\  (1-\min_{i\in L}\alpha_i(t))\sqrt{|L|}\max_{i\in L}\norm{x_i(0)-g}\to 0$$ as $t\to\infty.$ There is $t_0\geq 0$ such that $$\sqrt{2}\delta\min_{i\in L}\alpha_i(t)/|L|^4-(1-\min_{i\in L}\alpha_i(t))\sqrt{|L|}\max_{i\in L}\norm{x_i(0)-g}\geq \delta/|L|^4$$ for all $t\geq t_0.$ Assume that asymptotic stability does not hold in~\eqref{leader group model}. Then, there are $\delta>0$ and $(s_k)_{k\geq 0}$ increasing with $s_0\geq t_0$ and some component in profile $G(t_k)\cap\mathscr{G}(t_k)$ on~$L$ $\delta$-nontrivial for all $k\geq 0.$ Letting 
    $$M_0=4|L|^2\max_{i\in L}\norm{x_i(0)-g}\big(\max_{i\in L}\norm{x_i(0)-g}\vee\max_{i,j\in L}\norm{x_i(0)-x_j(0)}\big),$$ 
    it turns out from Lemma~\ref{lemma:Wt} that
    \begin{align*}
        W_0+M_0\sum_{t=0}^m(1/\min_{i\in L}\alpha_i(t)-1)&\geq\sum_{t=0}^m(W_t-W_{t+1})+M_0\sum_{t=0}^m(1/\min_{i\in L}\alpha_i(t)-1)\\
        &\geq 4\sum_{t=0}^m\sum_{i\in L}\norm{x_i(t)-x_i(t+1)}^2\ \hbox{for all}\ m\geq 0.
    \end{align*}
    As $m\to\infty,$ we derive
    \begin{align*}
        \infty&>W_0+M_0\sum_{t\geq 0} (1/\min_{i\in L}\alpha_i(t)-1)\geq 4\sum_{t\geq 0}\sum_{i\in L}\norm{x_i(t)-x_i(t+1)}^2\\
        &\geq 4\sum_{k\geq 0}\sum_{i\in L}\norm{x_i(s_k)-x_i(s_k+1)}^2\geq 4\sum_{k\geq 0}\delta^2/|L|^8=\infty,\ \hbox{a contradiction.}
    \end{align*}
    Hence, all components of profile $G\cap\mathscr{G}$ on~$L$ are $\delta$-trivial after some time for all $\delta>0.$

    Next, we claim that no components of profile $G\cap\mathscr{G}$ on~$L$ interact with each other after some time. It follows from claim~\ref{C1} that all components of profile $G\cap\mathscr{G}$ on~$L$ are $\epsilon/4$-trivial after some time~$s_0.$ Assume that claim~\ref{C2} is not the case. By finiteness of the social graph, there are edge $(i,j)$ and $(t_k)_{k\geq 0}$ increasing with $t_0\geq s_0$ such that vertices $i$ and $j$ belong to distinct components of profile $G\cap\mathscr{G}(t_k)$ on $L$, $$(i,j)\in E\cap \mathscr{E}(t_k)^c\ \hbox{and}\ (i,j)\in E\cap\mathscr{E}(t_k+1).$$ 
    Letting $y_i=y_i(t_k),\ y_i^\star=y_i(t_k+1)$ and $\alpha_i=\alpha_i(t_k)$ for all $i\in L$ and $k\geq 0$, it turns out from the triangle inequality that
    \begin{align*}
        \epsilon<\norm{y_i-y_j}&\leq \norm{y_i-y_i^\star/\alpha_i}+\norm{y_i^\star/\alpha_i-y_i^\star}+\norm{y_i^\star-y_j^\star}+\norm{y_j^\star-y_j^\star/\alpha_j}\\
        &\hspace{0.3cm}+\norm{y_j^\star/\alpha_j-y_j}\leq \epsilon/2+\norm{y_i^\star/\alpha_i-y_i^\star}+\norm{y_i^\star-y_j^\star}+\norm{y_j^\star-y_j^\star/\alpha_j}
    \end{align*}
    for the last inequality following from $\norm{y_i-y_i^\star/\alpha_i}\leq \epsilon/4$ and $\norm{y_j^\star/\alpha_j-y_j}\leq \epsilon/4.$
    Since $\limsup_{k\to\infty}\norm{y_i^\star/\alpha_i-y_i^\star}=0=\limsup_{k\to\infty}\norm{y_j^\star-y_j^\star/\alpha_j}$, we derive $$\epsilon/2 \leq \liminf_{k\to\infty}\norm{y_i^\star-y_j^\star}=\liminf_{k\to\infty}\norm{x_i^\star-x_j^\star},\ \hbox{a contradiction.}$$
\end{proof}

\eqref{leader group model} reduces to the synchronous Hegselmann-Krause model when \(\alpha_i(t) = 1\) for all \(i \in L\) at all times. Therefore, \(\sum_{t \geq 0}\left(1/\min_{i \in L} \alpha_i(t) - 1\right) = 0 < \infty\). From Theorem~\ref{Thm:asymptotic stability on L}, it follows that all components of a profile on $L$ become \(\epsilon\)-trivial, and no components interact with each other after some time under undirected social and opinion graphs. This indicates that all components achieve their consensus at the next time step, substantiating the finite time convergence property of the synchronous Hegselmann-Krause model under undirected social and opinion graphs.

\section{The follower group model}

Follower groups are independent. We first consider a leader group~$L$ and a follower group~$F$. \eqref{follower group model} becomes
$$x_i(t+1)=\frac{(1-\beta_i(t))}{|\N_i^F(t)|}\sum_{j\in \N_i^F(t)}x_j(t)+\frac{\beta_i(t)}{|\N_i^L(t)|}\sum_{j\in \N_i^L(t)}x_j(t),\quad i\in F,$$
which is equivalent to 
\begin{equation}\label{Follower model with a leader group}
    y_i(t+1)=\frac{(1-\beta_i(t))}{|\N_i^F(t)|}\sum_{j\in \N_i^F(t)}y_j(t)+\frac{\beta_i(t)}{|\N_i^L(t)|}\sum_{j\in \N_i^L(t)}y_j(t),\quad i\in F.
\end{equation}

\begin{lemma}\label{cs}
Assume that all followers have one social neighbor in leader group $L$ with target $g$, that $\big\{x_i(t)\big\}_{i\in L\cup  F}\subset \B(g,\min_{i\in F}\epsilon_i) $ for some $t\geq 0$ and that $$\sup_{s\geq t}\big\{\max_{i\in F}(1-\beta_i(s)),\ \max_{i\in L}\alpha_i(s)\big\}<1.$$ Then, $$\lim_{t\to \infty}\max_{i\in  L\cup F} \norm{x_i(t)-g}=0.$$
\end{lemma}

\begin{proof}
    $\big\{x_i(t)\big\}_{i\in L\cup  F}\subset \B(g,\min_{i\in F}\epsilon_i) $ for some $t\geq 0$ implies $\big\{x_i(s)\big\}_{i\in L\cup  F}\subset \B(g,\min_{i\in F}\epsilon_i) $ for all $s\geq t$. Via Theorem~\ref{Thm:consensus equal to target}, $\lim_{t\to\infty}\max_{k\in L}\norm{y_k(t)}=0$ therefore $$\max_{k\in L}\norm{y_k(s)}<\delta\ \hbox{for all $\delta>0$,  for some $p\geq t$}$$  and for all $s\geq p.$
For all $s\geq p$, $\delta>0$, $i\in F$ and $j\in L$, we have $$\norm{y_i(s)-y_j(s)}\leq \norm{y_i(s)}+\norm{y_j(s)}< \min_{i\in F}\epsilon_i+\delta $$ 
therefore $\norm{x_i(s)-x_j(s)}=\norm{y_i(s)-y_j(s)}\leq \min_{i\in F}\epsilon_i\leq \epsilon_i$ and $\N_i^L(s)\neq\emptyset$.
Let $\alpha_t=\max_{k\in L}\alpha_k(t)$, $\tilde{\beta}_t=\max_{k\in F}(1-\beta_k(t))$, $\gamma=\sup_{s\geq t}\{\tilde{\beta}_s,\ \alpha_s\}$,  $A_t=\max_{k\in F}\norm{y_k(t)}$ and $Z_t=\max_{k\in L}\norm{y_k(t)}$.
Applying the triangle inequality on~\eqref{Follower model with a leader group}, for all $i\in F$ and $t> p$,
\begin{align*}
  & A_{t+1}\leq\  \tilde{\beta}_t A_t+ Z_t\\
  &\hspace{0.4cm} \leq\  \tilde{\beta}_t\tilde{\beta}_{t-1}\ldots\tilde{\beta}_p A_{p}+\tilde{\beta}_t\ldots\tilde{\beta}_{p+1}Z_p+\ldots+\tilde{\beta}_t Z_{t-1} + Z_{t}\\
 &\hspace{0.4cm}  \leq\  \gamma^{t-p+1}A_p+(t-p+1)\gamma^{t-p}Z_p
\end{align*}
therefore
$$\limsup_{t\to\infty}A_{t+1}\leq 0.$$
This completes the proof.
\end{proof}

We move on to the follower group model with $m$ leader groups.

\begin{proof}[\bf Proof of Theorem~\ref{Thm:consensus on F}]
    $\big\{x_i(t),\ g_k\big\}_{i\in (\bigcup_{k=1}^m L_k)\cup F,\ k\in [m]}\subset \B(g_j,\min_{i\in F}\epsilon_i)$ for some $t\geq 0$ implies $$\big\{x_i(s)\big\}_{i\in (\bigcup_{k=1}^m L_k)\cup F}\subset \B(g_j,\min_{i\in F}\epsilon_i)\ \hbox{for all}\ s\geq t.$$ It follows from Theorem~\ref{Thm:consensus equal to target} that $\lim_{t\to\infty}\max_{k\in [m]}\max_{i\in L_k}\norm{x_i(t)-g_k}=0$ therefore $$\max_{k\in [m]}\max_{i\in L_k}\norm{x_i(s)-g_k}<\delta\ \hbox{for all}\ \delta>0,\  \hbox{for some}\ p\geq t\ \hbox{and for all}\ s\geq p.$$
For all $s\geq p$, $\delta>0$, $i\in F$ and $j\in L_k$, we have $$\norm{x_i(s)-x_j(s)}\leq \norm{x_i(s)-g_k}+\norm{g_k-x_j(s)}< \min_{i\in F}\epsilon_i+\delta $$ 
therefore $\norm{x_i(s)-x_j(s)}\leq \min_{i\in F}\epsilon_i\leq \epsilon_i$ and $\N_i^{L_k}(s)\neq \emptyset$. Letting
    $$\begin{array}{lll}
      \displaystyle\tilde{\beta}=\max_{i\in F}(1-\sum_{k=1}^m\beta_i^k),   & \displaystyle\gamma=\sup_{s\geq t}\big\{\max_{k\in [m]} \max_{i\in L_k}\alpha_i^k(s),\ \tilde{\beta}\big\},&\n \displaystyle g=\sum_{k=1}^m\beta_i^k g_k/\sum_{k=1}^m\beta_i^k,\vspace{2pt} \\
         \displaystyle A_t=\max_{i\in F}\norm{x_i(t)-g}/m, &\displaystyle C_t=\sum_{k=1}^m\max_{i\in L_k}\norm{x_i(t)-g_k}/m.
    \end{array}$$
Letting $$\bar{x}_i^F(t)=\frac{1}{|\N_i^F(t)|}\sum_{j\in \N_i^F(t)}x_j(t)\quad \hbox{and}\quad \bar{x}_i^{L_k}(t)=\frac{1}{|\N_i^{L_k}(t)|}\sum_{j\in \N_i^{L_k}(t)}x_j(t),$$
write~\eqref{follower group model} as
$$x_i(t+1)-g=(1-\sum_{k\in [m]}\beta_i^k)(\bar{x}_i^F(t)-g)+\sum_{k\in [m]}\beta_i^k(\bar{x}_i^{L_k}(t)-g_k)$$
and apply the triangle inequality, for all $i\in F$ and $t> p$,
\begin{align*}
  & A_{t+1}\leq\  \tilde{\beta}_t A_t+ C_t\\
  &\hspace{0.4cm} \leq\  \tilde{\beta}_t\tilde{\beta}_{t-1}\ldots\tilde{\beta}_p A_{p}+\tilde{\beta}_t\ldots\tilde{\beta}_{p+1}C_p+\ldots+\tilde{\beta}_t C_{t-1} + C_{t}\\
 &\hspace{0.4cm}  \leq\  \gamma^{t-p+1}A_p+(t-p+1)\gamma^{t-p}C_p
\end{align*}
therefore
$$\limsup_{t\to\infty}A_{t+1}\leq 0.$$
This completes the proof.
\end{proof}

\begin{lemma}\label{lemma:Xt}
    Assume that the social graph is undirected on $F$ and the opinion graph is undirected on $F$ with confidence threshold $\epsilon$. Let $X_t=\sum_{i,j\in F}(\norm{x_i(t)-x_j(t)}^2\wedge\epsilon^2)\vee\epsilon^2\mathbbm{1}\{(i,j)\notin E(t)\}$ and $E(t)\subset E(t+1).$ Then,
    \begin{align*}
        &X_t-X_{t+1}\geq 4\sum_{i\in F}\norm{x_i(t)-x_i(t+1)}^2- 4m|F|^2\max_{i\in F;k\in [m]}\beta_i^k(t)\\
        &\hspace{2cm}\times\bigg(\max_{i,j\in \bigcup_{k\in [m]}L_k\cup F}\norm{x_i(0)-x_j(0)}\vee\max_{i\in \bigcup_{k\in [m]}L_k\cup F}\norm{x_i(0)-g}\bigg)^2.
    \end{align*}
\end{lemma}

\begin{proof}
    Let $x_i=x_i(t)$, $x_i^\star=x_i(t+1)$, $\beta_i^k=\beta_i^k(t)$, $\N_i^F=\N_i^F(t)$, $\N_i^L=\N_i^L(t)$, $\Bar{x}_i^F=\sum_{j\in \N_i^F}x_j/|\N_i^F|$ and $\Bar{x}_i^L=\sum_{j\in \N_i^L}x_j/|\N_i^L|$. Observe that
    \begin{align*}
        &\hspace{-2pt} X_t-X_{t+1}\geq\sum_{i\in F}\sum_{j\in \N_i^F}(\norm{x_i-x_j}^2-\norm{x_i^\star-x_j^\star}^2)\\
        &=2\sum_{i\in F}\sum_{j\in \N_i^F}(<x_i-x_i^\star,x_i^\star-x_j>-<x_i^\star-x_j,x_j-x_j^\star>)\\
        &=2\sum_{i\in F}|\N_i^F|<x_i-x_i^\star,x_i^\star-\Bar{x}_i^F>+2\sum_{i\in F}\sum_{j\in \N_i^F}<x_i^\star-x_i,x_j^\star-x_j>\\
        &\hspace{1cm}-2\sum_{i\in F}\sum_{j\in \N_i^F}<x_i-x_j,x_j-x_j^\star>\\
        &\geq 2\sum_{i\in F}|\N_i^F|<x_i-x_i^\star,x_i^\star-\Bar{x}_i^F>+2\sum_{i\in F}\norm{x_i^\star-x_i}^2\\
        &\hspace{0.5cm}-\sum_{i\in F}\sum_{j\in \N_i^F-\{i\}}(\norm{x_i^\star-x_i}^2+\norm{x_j^\star-x_j}^2)-2\sum_{j\in F}|\N_j^F|<\Bar{x}_j^F-x_j,x_j-x_j^\star>\\
        &=4\sum_{i\in F}|\N_i^F|<x_i-x_i^\star,x_i^\star-\Bar{x}_i^F>+2\sum_{i\in F}\norm{x_i^\star-x_i}^2\\
        &\hspace{1cm}-2\sum_{i\in F}(|\N_i^F|-1)\norm{x_i-x_i^\star}^2+2\sum_{i\in F}|\N_i^F|\norm{x_i-x_i^\star}^2\\
        &=4\sum_{i\in F}\sum_{k\in [m]}\beta_i^k|\N_i^F|<x_i-x_i^\star,\Bar{x}_i^{L_k}-\Bar{x}_i^F>+4\sum_{i\in F}\norm{x_i-x_i^\star}^2\\
        &\geq 4\sum_{i\in F}\norm{x_i-x_i^\star}^2- 4m|F|^2\max_{i\in F;k\in [m]}\beta_i^k\\
        &\hspace{2cm}\times\bigg(\max_{i,j\in \bigcup_{k\in [m]}L_k\cup F}\norm{x_i(0)-x_j(0)}\vee\max_{i\in \bigcup_{k\in [m]}L_k\cup F}\norm{x_i(0)-g}\bigg)^2.
    \end{align*}
\end{proof}

\begin{lemma}\label{lemma: delta-nontrivial component on F}
    Assume that the social graph and opinion graph are undirected on $F$. If some component $H$ of the profile $G(t)\cap\mathscr{G}(t)$ on $F$ is $\delta$-nontrivial, then
    \begin{align*}
        \sqrt{\sum_{i\in F}\norm{x_i(t)-x_i(t+1)}^2}\geq& \sqrt{2}\delta(1-\max_{i\in F}\sum_{k\in [m]}\beta_i^k(t))/|F|^4-2m|F|\max_{k\in [m];i\in F}\beta_i^k(t)\\
        &\hspace{2cm}\times\left(\max_{i\in\bigcup_{k\in [m]}L_k\cup F}\norm{x_i(0)}\vee\max_{k\in [m]}\norm{g_k}\right).
    \end{align*}
\end{lemma}

\begin{proof}
    Letting \( V(H) \), the vertex set of \( H \), be \([h]\) and \(\mathbbm{1}=(1,\ldots,1)' \in \mathbb{R}^h\), express \(\mathbb{R}^h = W \oplus W^\perp\) for \( W = \Span(\{\mathbbm{1}\})\). For $x(t)=(x_1(t),\ldots,x_h(t))',$ write $$ x (t) = \left[c_1 \mathbbm{1} \,| \,c_2 \mathbbm{1} \,| \,\cdots \,| \,c_d \mathbbm{1} \right] +
            \left[\hat{c}_1 u^{(1)} \,| \,\hat{c}_2 u^{(2)} \,| \,\cdots \,| \,\hat{c}_d u^{(d)} \right] $$
    where $c_i$ and $\hat{c}_i$ are constants and $u^{(i)} \in \mathbbm{1}^\perp$ is a unit vector for all~$i \in [d]$. Observe that $$\norm{x_i(t)-x_j(t)}^2=\sum_{k\in [d]}\hat{c}_k^2(u^{(k)}_i-u^{(k)}_j)^2\leq 2\sum_{k\in [d]}\hat{c}_k^2\big((u^{(k)}_i)^2+(u^{(k)}_j)^2\big)\leq 2\sum_{k\in [d]}\hat{c}_k^2$$ for all $i,j\in [h].$ Hence, 
    $$\hbox{component}\ H\ \delta\hbox{-nontrivial implies}\ \sum_{k\in [d]}\hat{c}_k^2>\delta^2/2.$$ Letting $\tilde{\beta}(t)=(1-\sum_{k\in [m]}\beta_1^k(t),\ldots,1-\sum_{k\in [m]}\beta_h^k(t))'$ and $B(t)=\diag(\Tilde{\beta}(t))A(t)$ for $A(t)\in \mathbb{R}^{h\times h}$ with $A_{i,j}(t)=\mathbbm{1}\{j\in \N_i^F(t)\}/|\N_i^F(t)|,$ we get $$x(t)-x(t+1)=(I-B(t))x(t)-O(t)=\left[C(t)+F(t)\mathscr{L}(t)\right]x(t)-O(t)$$ where $C(t)=I-\diag(\Tilde{\beta}(t))$, $F(t)=\diag(\Tilde{\beta}(t))\left(\diag((d_i)_{i=1}^h)+I\right)^{-1}$ with $d_i$ the degree of vertex $i$ in component $H$, $O(t)=\sum_{k\in [m]}\diag((\beta_i^k(t))_{i\in [h]})(\bar{x}_i^{L_k})_{i\in [h]}'$ and $\mathscr{L}(t)$ is the Laplacian of component $H.$ It follows from Lemmas~\ref{L7} and~\ref{L10} that
    $$\lambda_2(\mathscr{L})>\frac{(2/h)^2}{2h}=2/h^3,$$
    \begin{align*}
        \norm{F(t)\mathscr{L}(t)x(t)}^2&=\sum_{k\in [d]}\hat{c}_k^2\norm{F(t)\mathscr{L}(t)u^{(k)}}^2\geq \sum_{k\in [d]}\hat{c}_k^2\lambda_2\left(\mathscr{L}(t)F^2(t)\mathscr{L}(t)\right)\\
        &\geq (\delta^2/2) (\min_{i\in [h]}\Tilde{\beta}_i(t)/h)^2\lambda^2_2(\mathscr{L}(t))\geq  2\delta^2(1-\max_{i\in [h]}\sum_{k\in [m]}\beta_i^k)^2/h^8.
    \end{align*}
    On the other hand, it follows from the triangle inequality that
    \begin{align*}
        \norm{C(t)x(t)}&\leq \sum_{i\in [h]}\sum_{k\in [m]}\norm{\beta_i^k(t)x_i(t)}\\
        &\leq mh\max_{i\in [h]; k\in [m]}\beta_i^k(t)\left(\max_{i\in \bigcup_{k\in [m]}L_k\cup F}\norm{x_i(0)}\vee \max_{k\in [m]}\norm{g_k}\right),\\
       \norm{O(t)} &\leq \sum_{i\in [h]}\sum_{k\in [m]}\norm{\beta_i^k\bar{x}_i^{L_k}}\\
       &\leq mh\max_{k\in [m];i\in [h]}\beta_i^k(t)\left(\max_{i\in\bigcup_{k\in [m]}L_k}\norm{x_i(0)}\vee\max_{k\in [m]}\norm{g_k}\right),
    \end{align*}
    therefore
    \begin{align*}
        \sqrt{\sum_{i\in F}\norm{x_i(t)-x_i(t+1)}^2}&\geq\sqrt{\sum_{i\in [h]}\norm{x_i(t)-x_i(t+1)}^2}=\norm{x(t)-x(t+1)}\\
        &\hspace{-3cm}=\norm{\left[F(t)\mathscr{L}(t)+C(t)\right]x(t)-O(t)}\geq \norm{F(t)\mathscr{L}(t)y(t)}-\norm{C(t)y(t)}-\norm{O(t)}\\
        &\hspace{-3cm}\geq \sqrt{2}\delta(1-\max_{i\in [h]}\sum_{k\in [m]}\beta_i^k)/h^4-2mh\max_{k\in [m];i\in [h]}\beta_i^k(t)\\
        &\hspace{2.5cm}\times\left(\max_{i\in\bigcup_{k\in [m]}L_k\cup F}\norm{x_i(0)}\vee\max_{k\in [m]}\norm{g_k}\right)\\
        &\hspace{-3cm}\geq \sqrt{2}\delta(1-\max_{i\in F}\sum_{k\in [m]}\beta_i^k)/|F|^4-2m|F|\max_{k\in [m];i\in F}\beta_i^k(t)\\
        &\hspace{2.5cm}\times\left(\max_{i\in\bigcup_{k\in [m]}L_k\cup F}\norm{x_i(0)}\vee\max_{k\in [m]}\norm{g_k}\right).
    \end{align*} 
\end{proof}


\begin{proof}[\bf Proof of Theorem~\ref{Thm:asymptotic stability}]
    We claim the following:
    \begin{enumerate}
        \item\label{C1F} All components of profile $G\cap\mathscr{G}$ on~$F$ are $\delta$-trivial after some time for all $\delta>0.$\vspace{2pt}
        \item\label{C2F} No components of profile $G\cap\mathscr{G}$ on~$F$ interact with each other after some time.
    \end{enumerate}
    Without loss of generality, we assume the social graph on $F$ remains constant over time, saying $G(t)|_F=G|_F=(F, E)$ for all $t\geq 0$. Observe that $$\sum_{t\geq 0}\max_{i\in F; k\in [m]}\beta_i^k(t)<\infty \implies \lim_{t\to\infty}\beta_i^k(t)=0\ \hbox{for all}\ i\in F\ \hbox{and}\ k\in [m].$$ Hence, we derive 
        \begin{align*}
            &a_t=\sqrt{2}\delta(1-\max_{i\in F}\sum_{k\in [m]}\beta_i^k(t))/|F|^4\to \sqrt{2}\delta/|F|^4,\\
            &b_t=2m|F|\max_{k\in [m];i\in F}\beta_i^k(t)\left(\max_{i\in\bigcup_{k\in [m]}L_k\cup F}\norm{x_i(0)}\vee\max_{k\in [m]}\norm{g_k}\right)\to 0\ \hbox{as}\ t\to\infty.
        \end{align*}
        There is $t_0\geq 0$ such that $$a_t-b_t\geq \delta/|F|^4\ \hbox{for all}\ t\geq t_0.$$  Assume that asymptotic stability does not hold in~\eqref{leader group model}. Then, there are $\delta>0$ and $(s_k)_{k\geq 0}$ increasing with $s_0\geq t_0$ and some component in profile $G(t_k)\cap\mathscr{G}(t_k)$ on~$F$ $\delta$-nontrivial for all $k\geq 0.$ Letting 
    $$M_0=4m|F|^2\bigg(\max_{i,j\in \bigcup_{k\in [m]}L_k\cup F}\norm{x_i(0)-x_j(0)}\vee\max_{i\in \bigcup_{k\in [m]}L_k\cup F}\norm{x_i(0)-g}\bigg)^2,$$ 
    it turns out from Lemma~\ref{lemma:Xt} that
    \begin{align*}
        X_0+M_0\sum_{t=0}^{\hat{m}}\max_{i\in F; k\in [m]}\beta_i^k(t)&\geq\sum_{t=0}^{\hat{m}}(X_t-X_{t+1})+M_0\sum_{t=0}^{\hat{m}}\max_{i\in F; k\in [m]}\beta_i^k(t)\\
        &\geq 4\sum_{t=0}^{\hat{m}}\sum_{i\in F}\norm{x_i(t)-x_i(t+1)}^2\ \hbox{for all}\ \hat{m}\geq 0.
    \end{align*}
    As $\hat{m}\to\infty,$ we derive
    \begin{align*}
        \infty&>W_0+M_0\sum_{t\geq 0} \max_{i\in F; k\in [m]}\beta_i^k(t)\geq 4\sum_{t\geq 0}\sum_{i\in F}\norm{x_i(t)-x_i(t+1)}^2\\
        &\geq 4\sum_{k\geq 0}\sum_{i\in F}\norm{x_i(s_k)-x_i(s_k+1)}^2\geq 4\sum_{k\geq 0}\delta^2/|F|^8=\infty,\ \hbox{a contradiction.}
    \end{align*}
    Hence, all components of profile $G\cap\mathscr{G}$ on~$F$ are $\delta$-trivial after some time for all $\delta>0.$

    Next, we claim that no components of profile $G\cap\mathscr{G}$ on~$F$ interact with each other after some time. It follows from claim~\ref{C1F} that all components of profile $G\cap\mathscr{G}$ on~$F$ are $\epsilon/4$-trivial after some time~$s_0.$ Assume that claim~\ref{C2F} is not the case. By finiteness of the social graph, there are edge $(i,j)$ and $(t_k)_{k\geq 0}$ increasing with $t_0\geq s_0$ such that vertices $i$ and $j$ belong to distinct components of profile $G\cap\mathscr{G}(t_k)$ on $F$, $$(i,j)\in E\cap \mathscr{E}(t_k)^c\ \hbox{and}\ (i,j)\in E\cap\mathscr{E}(t_k+1).$$ 
    Letting 
    $$\begin{array}{lll}
         \displaystyle \Bar{x}_i^F=\frac{1}{|\N_i^F(t_k)|}\sum_{j\in \N_i^F(t_k)}x_j(t_k),&\displaystyle x_i=x_i(t_k),&\displaystyle x_i^\star=x_i(t_k+1),\vspace{2pt}\\
         \displaystyle\Bar{x}_i^L=\frac{1}{|\N_i^L(t_k)|}\sum_{j\in \N_i^L(t_k)}x_j(t_k),&\displaystyle \Tilde{\beta}_i^j=\beta_i^j(t_k), &\displaystyle\Tilde{\beta}_i=1-\sum_{j\in [m]}\beta_i^j(t_k)
    \end{array}$$
    for all $i\in F$ and $k\geq 0$, it turns out from the triangle inequality that
    \begin{align*}
        \epsilon&<\norm{x_i-x_j}\leq \norm{x_i-x_i^\star}+\norm{x_i^\star-x_j^\star}+\norm{x_j^\star-x_j}.
    \end{align*}
    On top of that, we get
    \begin{align*}
        \norm{x_i-x_i^\star}&\leq \Tilde{\beta}_i\norm{x_i-\bar{x}_i^F}+\norm{\sum_{j\in [m]}\beta_i^j(x_i-\bar{x}_i^L)}\\
        &\leq \Tilde{\beta}_i\epsilon/4+m\max_{j\in [m]; i\in F}\beta_i^j\\
        &\hspace{0.5cm}\times\bigg(\max_{i,j\in \bigcup_{k\in [m]}L_k\cup F}\norm{x_i(0)-x_j(0)}\vee\max_{i\in \bigcup_{k\in [m]}L_k\cup F}\norm{x_i(0)-g}\bigg),
    \end{align*}
    similarly for $\norm{x_j-x_j^\star},$ therefore 
    $$\liminf_{k\to\infty}\norm{x_i-x_i^\star}\leq \epsilon/4\quad \hbox{and}\quad \liminf_{k\to\infty}\norm{x_j-x_j^\star}\leq \epsilon/4.$$
   This implies $$\epsilon/2 \leq \liminf_{k\to\infty}\norm{x_i^\star-x_j^\star},\ \hbox{a contradiction.}$$
   It follows from Claims~\ref{C1F} and~\ref{C2F} that $\sum_{j\in \N_i^F(t)}x_j(t)/|\N_j^F(t)|$ converges to some random variable $\Tilde{x}_i$ as $t\to\infty$ for all $i\in F.$ Since $\max_{k\in [m];i \in F}\beta_i^k(t)\to 0$ as $t\to\infty$ and $\sum_{j\in \N_i^L(t)}x_j(t)/|\N_i^L(t)|$ is bounded by $(\max_{i\in\bigcup_{k\in [m]}L_k}\norm{x_i(0)}\vee\max_{k\in [m]}\norm{g_k}),$ we get $x_i(t+1)\to \Tilde{x}_i$ as $t\to\infty$ for all $i\in F.$
\end{proof}

\section*{Acknowledgment}
The author is partially funded by NSTC grant.

\end{document}